\documentclass{article}
\usepackage{textcomp}
\usepackage{amsmath}
\usepackage{subfig}
\usepackage[colorlinks=true, pdfstartview=FitV, linkcolor=black, citecolor= black, urlcolor= black]{hyperref}
\usepackage{footnpag}	
\usepackage{amsfonts}
\usepackage{amssymb}
\usepackage{gensymb}

\usepackage{subfig}
\usepackage{subcaption}
\usepackage{graphicx} 
\usepackage{algorithm}
\usepackage{algpseudocode}
\usepackage{blindtext}
\usepackage{multicol}

\usepackage{arxiv}

\usepackage[utf8]{inputenc} 
\usepackage[T1]{fontenc}    
\usepackage{hyperref}       
\usepackage{url}            
\usepackage{booktabs}       
\usepackage{amsfonts}       
\usepackage{nicefrac}       
\usepackage{microtype}      
\usepackage{lipsum}		
\usepackage{doi}
\usepackage{caption}

\usepackage{rotating}
\usepackage[toc,page]{appendix}

\usepackage[numbers]{natbib}

\title{An Enhanced Shifted QR Algorithm for Efficient Eigenvalue Computation of Square Non-Hermitian Matrices}

\author{
\href{https://orcid.org/0000-0000-0000-0000}{\includegraphics[scale=0.06]{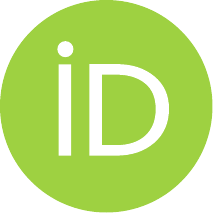}\hspace{1mm}Chahat Ahuja}\thanks{B.Tech CSAM, Department of Mathematics, IIIT Delhi} \\
\texttt{chahat22138@iiitd.ac.in} \\
\And
\href{https://orcid.org/0000-0000-0000-0000}{\includegraphics[scale=0.06]{orcid.pdf}\hspace{1mm}Partha Chowdhury}\thanks{Ph.D. Student, IIIT Delhi} \\
\texttt{parthac@iiitd.ac.in} \\
\And
\href{https://orcid.org/0000-0000-0000-0000}{\includegraphics[scale=0.06]{orcid.pdf}\hspace{1mm}Subhashree Mohapatra}\thanks{Assistant Professor, Department of Mathematics, IIIT Delhi} \\
\texttt{subhashree@iiitd.ac.in} \\
}

\date{}

\renewcommand{\headeright}{}
\renewcommand{\undertitle}{}

\hypersetup{
pdftitle={A template for the arxiv style},
pdfsubject={q-bio.NC, q-bio.QM},
pdfauthor={David S.~Hippocampus, Elias D.~Striatum},
pdfkeywords={First keyword, Second keyword, More},
}

\begin{document}
\maketitle

\begin{abstract}
	This work presents a novel approach to compute the eigenvalues of non-Hermitian matrices using an enhanced shifted QR algorithm. The existing QR algorithms fail to converge early in the case of non-hermitian matrices, and our approach shows significant improvement in convergence rate while maintaining accuracy for all test cases. In this work, though our prior focus will be to address the results for a class mid- large sized non-Hermitian matrices, our algorithm has also produced significant improvements in the case of comparatively larger matrices such as $50 \times 50$ non-Hermitian matrices.

\vspace{0.2cm}
\end{abstract}


\section{Introduction}
The eigenvalue problem involves finding a scalar $\lambda$ for a corresponding vector $\mathbf{v}$ such that $A\mathbf{v} = \lambda \mathbf{v}$ holds good. But when matrices become extremely large or ill-conditioned then the eigen value problem faces difficulties to reach its desired results. QR factorization being one of the methodologies to solve an eigen value problem, invovles in finding a orthogonal matrix $Q$ and an upper triangular matrix $R$ for a given matrix $A$(\cite{watkins_qr-like_2000},\cite{bini_fast_2005}). 
For Hermitian matrices where $A= A^H$ ( $A^H$ is the Conjugate transpose of A), efficient QR algorithms exist because of their well-structured spectral properties, including real eigenvalues and orthogonal eigenvectors \cite{nikpour_novel_2003}.
However, for non-Hermitian matrices, the problem becomes more challenging because of its complex eigenvalues and non-orthogonal eigenvectors. In this paper we have proposed an algorithm to address the drawbacks of QR algorithm in case of non-Hermitian matrices. Our algorithm performs significantly well in comparison with existing QR algorithms. While existing algorithms fail to achieve satisfactory results due to their slower convergence and less accuracy, the proposed algorithm addresses these challenges, making the convergence faster and maintaining the accuracy of eigenvalues for randomly generated non-Hermitian matrices.

\section{Problem Statement}
Given a $n \times n$  non-Hermitian  matrix $A \in\mathbb{C}^{n \times n}$, the goal is to compute the eigenvalues $\lambda_1, \lambda_2, \lambda_3 , \lambda_4 ....\in \mathbb{C}$ such that, 
\begin{equation*}
    A v_i = \lambda_i v_i \ \forall \ i=1,2,3,4,5,.....,n
\end{equation*}
where $v_i \in \mathbb{C}^n$ are the corresponding eigenvectors.
Our goal is to develop a fast and efficient algorithm for computing the eigenvalues of non-Hermitian matrices while ensuring numerical stability. The proposed method aims to reduce computational complexity by accelerating convergence using advanced shift strategies and early deflation techniques. By efficiently handling large and complex matrices, our approach minimizes iteration count and improves robustness against numerical errors. These enhancements make the algorithm well-suited for practical applications requiring high-speed and accurate eigenvalue computations.

\section{Brief Literature Survey}
 The QR algorithm remains one of the most fundamental techniques for computing eigenvalues of matrices, and various refinements have been proposed over the years. The classical QR method for determining all eigenvalues of real square matrices was extensively studied by Francis and Kublanovskaya, whose work laid the foundation for modern implementations \cite{fekadie_anley_qr_2016}. The convergence properties of the QR method, particularly in the context of shifted iterations for varied dimension matrices, have been analysed using structured matrices. Batterson\cite{batterson_convergence_nodate} examined the behaviour of the shifted QR algorithm for 3 × 3 standard matrices, highlighting cases where convergence issues arise. The restarted QR algorithm proposed by Calvetti \cite{calvetti_restarted_2002} demonstrated significant efficiency improvements for structured matrices by periodically restarting iterations to enhance numerical stability.
Meanwhile, Braman et.al. \cite{braman_multishift_2002} introduced aggressive early deflation techniques for multi-shift QR algorithms, significantly reducing computational cost by identifying and deflating converged eigenvalues early. The convergence characteristics of multi-shift QR algorithms have also been examined for symmetric tridiagonal matrices, revealing optimal shift selection strategies to improve performance \cite{su_convergence_2013}. Additionally, Ahu and Tisseur proposed a novel deflation criterion to enhance the robustness of multi-shift QR iterations\cite{ahues_new_1997}. Parlett \cite{parlett_necessary_1966} derived necessary and sufficient conditions for converging the QR algorithm on Hessenberg matrices, contributing to a deeper theoretical understanding of QR iterations. These studies form the basis for developing our Enhanced shifted QR algorithm, which aims to achieve faster convergence for non-Hermitian matrices while addressing the limitations of existing literature.

\section{Background}
The classical QR algorithm is an iterative method for computing the eigenvalues. It involves factorizing \( A \) into a product of an orthogonal matrix \( Q \) and an upper triangular matrix \( R \),
    Then, the matrix \( A \) is updated as:
    \[A_{k} =  Q R\]
    \[A_{k+1} = R Q\]
    The algorithm iterates this process until the matrix \( A_k \) converges to a diagonal matrix whose diagonal elements are the eigenvalues of \( A \).QR Decomposition can be performed using three methods: Householder QR Factorization (uses Householder reflections for numerical stability), Gram-Schmidt QR Factorization (constructs an orthonormal basis, available in classical and modified forms), and Givens Rotation QR Factorization (uses Givens rotations, suitable for sparse matrices).

\subsection{Shifted QR Decomposition}
The shifted QR algorithm \cite{robles_implementing_nodate} enhances the classic QR method. It accelerates the  algorithm's convergence by applying a shift to the matrix \( A_k \) at each iteration. The idea is to use  a shift \( \sigma_k \) that improves the convergence of the eigenvalues.
    \[A_k = Q_k R_k\]
    \[A_{k+1} = R_k Q_k + \sigma_k I\]
where \( \sigma_k \) is a scalar shift and \( I \) is the identity matrix. Basically, shifts accelerate convergence to the dominant eigenvalues reducing the off-diagonal entries faster. 

Now we will define the shifting strategy i.e the wilkinson shift that have been used in the algorithm.\\ 
Wilkinson Shift\cite{robles_implementing_nodate}: Let $B$ be the lower-right $2 \times 2$ submatrix of $A^{(k)}$. The Wilkinson shift is the eigenvalue of $B$ closer to $a_m$, defined as:
    \[
    \mu = a_m - \operatorname{sign} (\delta) b_{m-1}^2 \left( |\delta| + \sqrt{\delta^2 + b_{m-1}^2} \right),
    \]
    where\[
\delta = \frac{a_{m-1}}{2} - a_m
\] and $\operatorname{sign}(\delta)$ can be chosen arbitrarily if $\delta = 0$. This shift accelerates the QR algorithm to cubic convergence.

\subsubsection{Deflation} 
Deflation \cite{robles_implementing_nodate} is a technique used to reduce the size of a matrix, thereby decreasing computational complexity and improving efficiency. It eliminates already converged eigenvalues, progressively reducing the matrix size $m$. Here, we impose deflation on the off-diagonal elements so they can be deflated correctly in a certain tolerance threshold, making the algorithm computationally efficient.

Suppose, after several iterations, a matrix becomes the following form:
    \begin{center}
        $\begin{bmatrix}
   \lambda_1 & *         & * \\
   0         & \lambda_2 & * \\
   0         &  0        & \lambda_3\\
        \end{bmatrix}$
    \end{center}
where the sub-diagonal elements are effectively zero $(\approx 10 ^ {-12})$ [Ex:$A_{32}= 1.2 * 10^{-12}$]

At this point, the subdiagonal norm would be $\approx 0$ which is less than $10^{-10}$. The matrix is now nearly upper triangular, so we can deflate extracting $\lambda_3$ as an approximate eigenvalue. We also reduced the matrix size by removing the last row and column.
\section {Methodology}
The Enhanced Shifted QR Algorithm \footnote[1]{Python code of the algorithm available in \href{https://github.com/chahat6606/Enhanced-Shifted-QR-Algorithm}{\textbf{https://github.com/chahat6606/Enhanced-Shifted-QR-Algorithm}}} computes the eigenvalues of a matrix by iteratively applying the QR decomposition with an implicit Wilkinson shift for accelerated convergence. Additionally, after every iteration, it incorporates deflation by eliminating converged eigenvalues, progressively reducing the matrix size.

\subsection{Algorithm}
Firstly, the algorithm incorporates the Wilkinson shift, a well-known shift strategy focusing on the eigenvalues of $2 \times 2$ sub-matrix to ensure faster convergence towards dominant eigenvalues. We also employed early deflation techniques to accelerate computation, reducing the adequate size of the matrix by removing converged eigenvalues from the matrix. We have also incorporated a balancing step into the algorithm to ensure the algorithm's numerical stability, minimizing sensitivity to rounding errors during iterations. These enhancements make our algorithm novel compared to the existing ones, overcoming the traditional limitations. The algorithm first aims to converge in fewer iterations than existing algorithms to its dominant eigenvalues for matrices with complex spectra. Secondly, the algorithm reduces computational complexity by efficiently deflating converged eigenvalues. Lastly, The algorithm maintains numerical stability throughout the iterative process (relative errors). Here we are providing the pseudocode and the details algorithm is given in the  appendix \ref{appendix:Algorithm1}.
\subsection{Pseudocode of the Algorithm}
\textbf{Input:} 
\begin{itemize}
    \item Matrix \( A \in \mathbb{R}^{n \times n} \)
    \item Maximum iterations \( k_{\max} \)
    \item Convergence tolerance \( \epsilon \)
    \item Deflation tolerance \( \delta \)
\end{itemize}
\textbf{Output:} 
\begin{itemize}
    \item Eigenvalues of \( A \)
    \item Iterations required for convergence
    \item Off-diagonal norms
    \item Deflation count
\end{itemize}

\subsection{Time and Space Complexity}
In this Subsection we will discuss about the time and space complexity of the Enhanced Shifted QR algorithm in details.

The enhanced shifted QR algorithm is an iterative method for computing the eigenvalues of a given matrix $A$. 
In each iteration, the most computationally intensive operations involve orthogonal decomposition, which, when implemented using standard Gram-Schmidt or Householder reflections, requires $O(n^3)$ time for an $n \times n$ matrix. 
Since convergence usually requires $O(n)$ iterations in practice, the worst-case time complexity of the full algorithm is $O(n^4)$. 
However, empirical observations suggest that convergence often occurs in less than the worst-case bound, achieving at least $O(n^3)$ complexity in practice.

The algorithm operates in place on the matrix \( A \), iteratively modifying it. The additional space required is for the matrices \( Q \) and \( R \), both of which are \( O(n^2) \) in size. Thus, the overall space complexity is $O(n^2)$. Since no additional extensive data is stored, the space complexity remains manageable, making the algorithm feasible for moderately large matrices.

\subsection{Performance Analysis on matrices of arbitrary dimension}
We carried out extensive experiments using our algorithm with square non-hermitian matrices of various dimensions and summarized the findings as follows. The performance criteria were assessed with regard to the convergence speed, accuracy, and computational efficiency, tested across diverse matrix sizes.
\subsubsection{\textbf{Case 1: For \texorpdfstring{$3 \times 3$}{3x3} Matrices}}
\begin{figure}[H]
    \centering
    \includegraphics[width=1.0\linewidth]{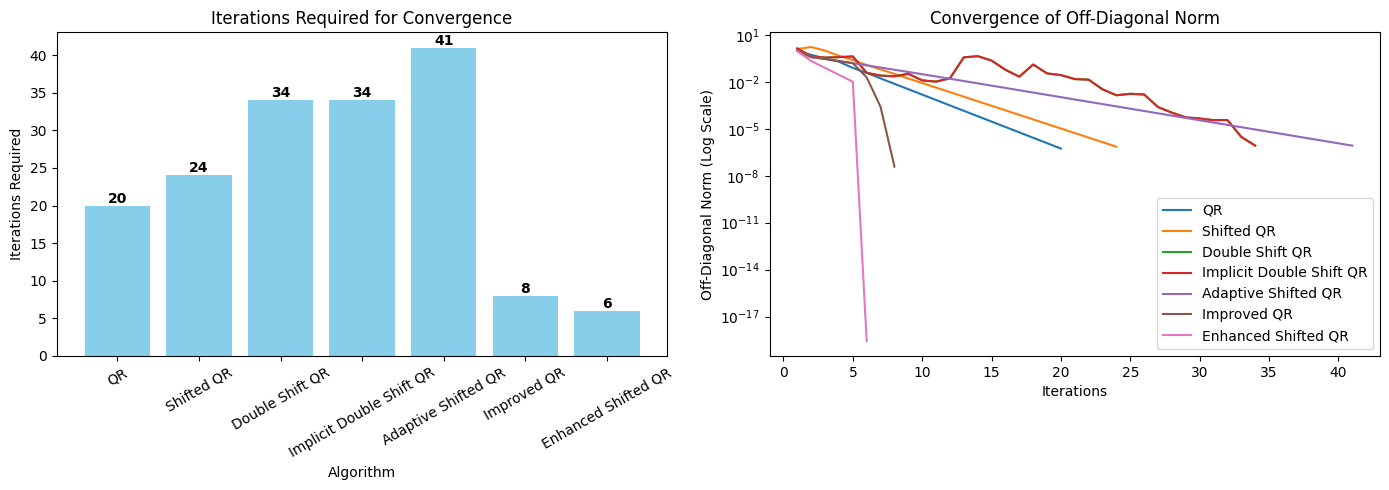}
    \begin{minipage}{0.9\linewidth}
        \centering
        (a) Iterations Required for Convergence \hspace{0.1cm} (b) Convergence of offdiagonal norm
    \end{minipage}
    
    \caption{Iterations Required and Off-Diagonal Norm Convergence of a $3 \times 3$ randomly generated non-Hermitian matrix.}
    \label{fig:q1_final.png}
\end{figure}

Firstly, we tested the Enhanced Shifted QR algorithm on various randomly generated $3 \times 3$ non-Hermitian matrices,
comparing its performance with all other existing algorithms. While all test cases followed a general trend, we present a
representative example in Figure ~\ref{fig:q1_final.png}, showcasing a randomly generated matrix and its results.

Figure ~\ref{fig:q1_final.png}(a) presents a comparative analysis of different QR-based methods in terms of convergence rate and off-diagonal norm reduction. It illustrates the number of iterations required for convergence, revealing that the Enhanced Shifted QR method exhibits the most rapid convergence, requiring only six iterations. In contrast to other algorithms which like Implicit Double Shift QR and Shifted QR methods that demand significantly higher iteration counts, 41 and 24, respectively, indicating slower convergence.

Figure~\ref{fig:q1_final.png}(b) demonstrates the convergence behavior of the off-diagonal norm across various iterations on a logarithmic scale. The Enhanced Shifted QR method shows a steep decline, suggesting a more efficient reduction of the off-diagonal norm compared to other methods, meaning that the matrix convergences to a diagonal form quickly.

Overall, our algorithm demonstrates the fastest performance for randomly generated $3 \times 3$ square non-Hermitian matrices. While other algorithms exist with nearly comparable efficiency in this low-dimensional setting, the difference in computational speed becomes significantly more pronounced as the matrix dimension increases, making our approach particularly advantageous for high-dimensional eigenvalue computations.

\subsubsection{\textbf{Case 2: For \texorpdfstring{\( 7 \times 7 \)}{7 x 7} Matrices}}
We gradually increased the dimension of the matrices to test the performance of the algorithm. Further,we check the performance of our algorithm for randomly generated $7 \times 7$ matrices to observe emerging trends. The following is a representation of its performance in comparison with the existing ones.
\begin{figure}[H]
    \centering
    \includegraphics[width=1.0\linewidth]{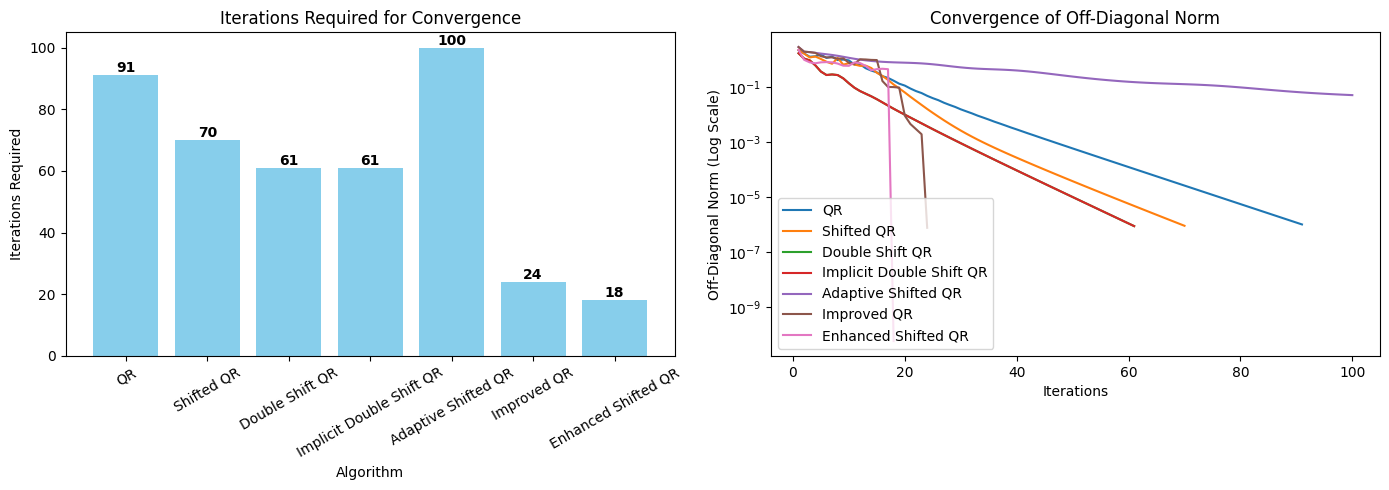}
    
    \begin{minipage}{0.9\linewidth}
        \centering
        (a) Iterations Required for Convergence \hspace{2cm} (b) Convergence of Off-Diagonal Norm
    \end{minipage}
    
    \caption{Iterations Required and Off-Diagonal Norm Convergence of a $7 \times 7$ randomly generated non-Hermitian matrix.}
    \label{fig:q2_final.png}
\end{figure}
 Figure~\ref{fig:q2_final.png}(a) shows the number of iterations required for convergence, indicating that our proposed method (Enhanced Shifted QR) requires 18 iterations, which is significantly less in comparison to standard QR methods requiring around 50-100 iterations. The Improved-QR algorithm demonstrates comparable performance to our method; however, as the dimensionality increases, the performance gap becomes more pronounced

Figure~\ref{fig:q2_final.png}(b) further validates this observation by illustrating the convergence rate of the off-diagonal norm. Our approach demonstrates the steepest decline, confirming the fast rate of the matrix of transforming into a diagnol form.

Finally, Our algorithm shows substantial improvement over traditional approaches, and its efficiency gain scales favorably with increasing matrix size. This makes it a highly promising method for high-dimensional eigenvalue computations, where computational cost becomes a critical factor.

\subsubsection{\textbf{Case 1: For \texorpdfstring{$50 \times 50$}{50x50} Matrices}}
To obtain a more general understanding of its scalability and effectiveness, we directly test it on a significantly larger matrix of size $50 \times 50$. This approach allows us to observe the overall trend rather than focusing on case-specific behavior. The results from this high-dimensional setting provide a clearer insight into the algorithm's performance improvements as matrix size increases.
\begin{figure}[H]
    \centering
    \includegraphics[width=1.0\linewidth]{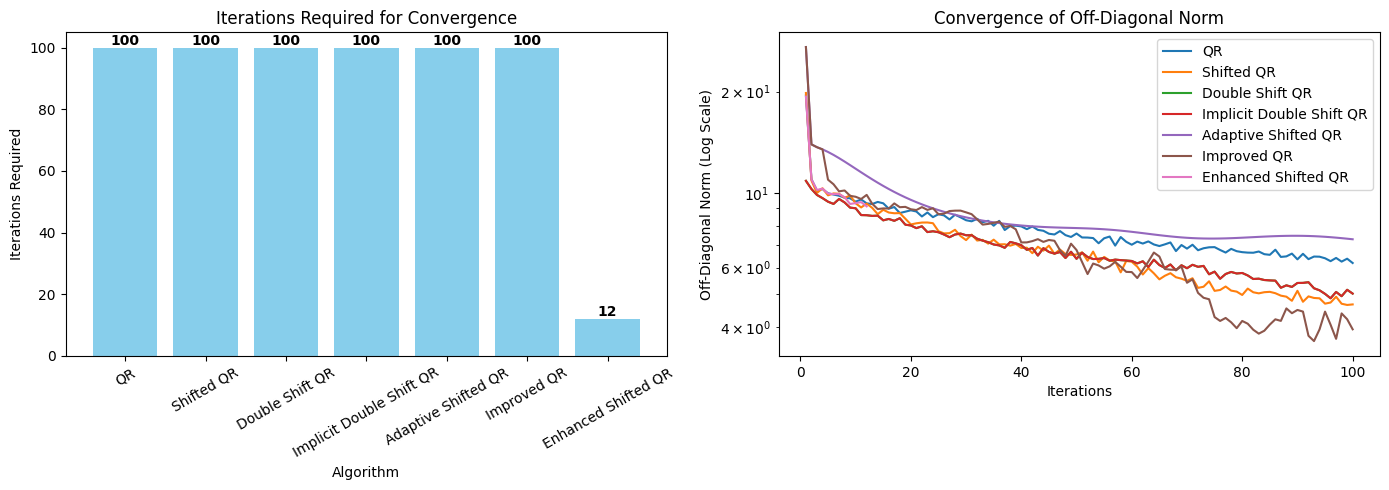}
    
    \begin{minipage}{0.9\linewidth}
        \centering
        (a) Iterations Required for Convergence \hspace{2cm} (b) Convergence of Off-Diagonal Norm
    \end{minipage}
    
    \caption{Iterations Required and Off-Diagonal Norm Convergence of a $50 \times 50$ randomly generated non-Hermitian matrix.}
    \label{fig:q3_final.png}
\end{figure}

Figure~\ref{fig:q3_final.png}(a) illustrates the performance of various QR-based algorithms for a larger matrix. The left plot clearly shows that while most traditional QR-based methods require the maximum iteration count to converge, our proposed Enhanced Shifted QR method converges in exceptionally fewer iterations, demonstrating its superior efficiency. 

Figure~\ref{fig:q3_final.png}(b) further reinforces this by depicting the rapid decay of the off-diagonal norm for our method, indicating faster convergence. This highlights that as matrix size increases, our approach not only outperforms others but also maintains a substantial efficiency advantage, making it a robust choice for high-dimensional eigenvalue computations.

Overall, the proposed algorithm consistently demonstrates a superior convergence rate, with the performance gap widening as the matrix dimension increases. Notably, it achieves the same level of accuracy as existing methods while significantly reducing the number of iterations required for convergence. This efficiency gain becomes increasingly pronounced in higher-dimensional settings, establishing our approach as a robust and scalable solution for eigenvalue computations.

\section{Results and Findings}
The results clearly highlight the efficiency of the Enhanced Shifted QR algorithm over traditional QR-based methods, particularly as the matrix size increases. For smaller matrices, such as \(3 \times 3\), our method achieves convergence in just \textbf{6 iterations}, while the Adaptive Shifted QR and Shifted QR methods require \textbf{41} and \textbf{24} iterations, respectively. This rapid convergence is also reflected in the off-diagonal norm reduction, where our approach consistently exhibits a steeper decline. As the matrix dimension increases to \(7 \times 7\), the efficiency gap widens further, with our method converging in \textbf{18 iterations}, compared to the \textbf{70-100} iterations observed for traditional QR methods. The Improved QR method performs relatively well, but our approach maintains a distinct advantage in both iteration count and norm reduction. When tested on large-scale \(50 \times 50\) matrices, the Enhanced Shifted QR method continues to outperform others by requiring substantially fewer iterations, often well below the \textbf{100+} iterations observed in standard QR-based techniques. The rapid decay of the off-diagonal norm in the log-scale plots further validates the robustness of our approach, demonstrating that it efficiently drives the matrix towards a diagonal form with minimal computational effort. The superior performance of our algorithm becomes increasingly pronounced as matrix size grows, making it an ideal choice for high-dimensional eigenvalue problems. Overall, our method not only ensures faster convergence but also maintains numerical stability, offering a significant computational advantage over existing algorithms.

\section{Conclusion}
The Enhanced Shifted QR algorithm establishes itself as a superior alternative to traditional QR-based eigenvalue solvers by consistently demonstrating rapid convergence and robust numerical stability. The significant reduction in iteration count across different matrix dimensions confirms the efficiency of our method, particularly in high-dimensional cases where computational complexity is a major concern. Traditional QR and Shifted QR methods exhibit a slower rate of off-diagonal norm decay, often requiring a substantially higher number of iterations to achieve a comparable diagonal form. The superior performance of our algorithm in both iteration reduction and norm minimization makes it an attractive choice for large-scale problems, where efficiency is paramount.

The future scope of this algorithm extends towards optimizing its implementation for parallel computing environments, allowing for even faster execution in high-performance computing systems. Additionally, incorporating adaptive heuristics for dynamic shift selection could further refine the convergence behavior, making the method even more robust for ill-conditioned matrices. Potential applications include quantum computing, spectral clustering, and numerical solutions to large-scale differential equations, where rapid eigenvalue decomposition is essential. Furthermore, extending the method to non-square and structured matrices could unlock new possibilities in data science and applied mathematics.

In conclusion, our Enhanced Shifted QR algorithm not only outperforms existing methods in terms of convergence speed and accuracy but also demonstrates scalability and robustness, making it a highly effective tool for modern eigenvalue computations. As computational demands continue to grow, our approach provides a promising direction for the development of even more efficient and powerful numerical techniques.




\begin{appendices}
\section{Appendix 1: Enhanced Shifted QR Algorithm for Eigenvalue Computation }

\begin{algorithm}
    \begin{algorithmic}[1]
        \Function{EnhancedShiftedQR}{$A, k_{\max}, \epsilon, \delta$}\label{appendix:Algorithm1}
        \State $A_k \gets A$
        \State $n \gets \text{size}(A)$
        \State {Initialize an empty list for eigenvalues}
        
        \For{$i = 1$ to $k_{\max}$}
            \For{$j = n-1$ down to $1$}
                \If{$\|A_k[j, :j]\| < \delta$} \Comment{Small sub-diagonal element}
                    \State Append $A_k[j, j]$ to eigenvalues
                    \State Remove the $j$th row and column from $A_k$
                    \State Update $n \gets n-1$
                    \State \textbf{Break} \Comment{Restart iteration with reduced size}
                \EndIf
            \EndFor

            \If{$n == 1$} 
                \State Append the last remaining diagonal element $A_k[0,0]$ to eigenvalues
                \State \textbf{Break}
            \EndIf

            \State Compute Wilkinson shift $\mu$ from the last $2 \times 2$ submatrix
            \State Compute QR decomposition: $Q, R \gets \text{QR decomposition of } (A_k - \mu I)$
            \State Update matrix: $A_k \gets RQ + \mu I$

            \If{$\|A_k\|_{\text{off-diagonal}} < \epsilon$} 
                \State Append all diagonal elements of $A_k$ to eigenvalues
                \State \textbf{Break}
            \EndIf
        \EndFor

        \State \Return Eigenvalues, iteration count, deflation count
        \EndFunction
    \end{algorithmic}
\end{algorithm}
\end{appendices}







\end{document}